\documentclass{article}
\usepackage[utf8]{inputenc}
\usepackage[a4paper, total={6.4in, 9.7in},top=1in]{geometry}
\usepackage{amsmath}
\usepackage{amsfonts}
\usepackage{amssymb}
\usepackage{amsthm}
\usepackage{hyperref}
\usepackage{parskip}
\usepackage{xcolor}
\usepackage[pdftex]{graphicx}
\usepackage{enumerate}
\usepackage{flafter}
\usepackage{listings}
\usepackage[backend=biber, style=authoryear, natbib=true, maxbibnames=4, maxcitenames=2]{biblatex}
\DeclareNameAlias{sortname}{last-first}

\DefineBibliographyStrings{english}{in = {{}{}}}

\newcommand{\mytitle}{Bayesian Fitting of Dirichlet Type I and II Distributions}
\hypersetup{pdfinfo={
	Title={\mytitle},
	Subject={\mytitle},
	Author={Sean van der Merwe},
	Keywords={Bayes, Beta, Dirichlet, Dirichlet Type II, Objective prior, Simulation} } }

\newlength{\figwidth} 

\usepackage{enumerate}

\usepackage{authblk}
\title{Bayesian Fitting of Dirichlet Type I and II Distributions}
\author[a,*]{Sean~van~der~Merwe}
\author[a]{Daan~de~Waal}
\affil[a]{University of the Free State, Box 339, Bloemfontein, 9300, South Africa}
\affil[*]{Corresponding author --- vandermerwes@ufs.ac.za}

\begin{document}
\maketitle

\begin{center}
\includegraphics[width=1.2\figwidth]{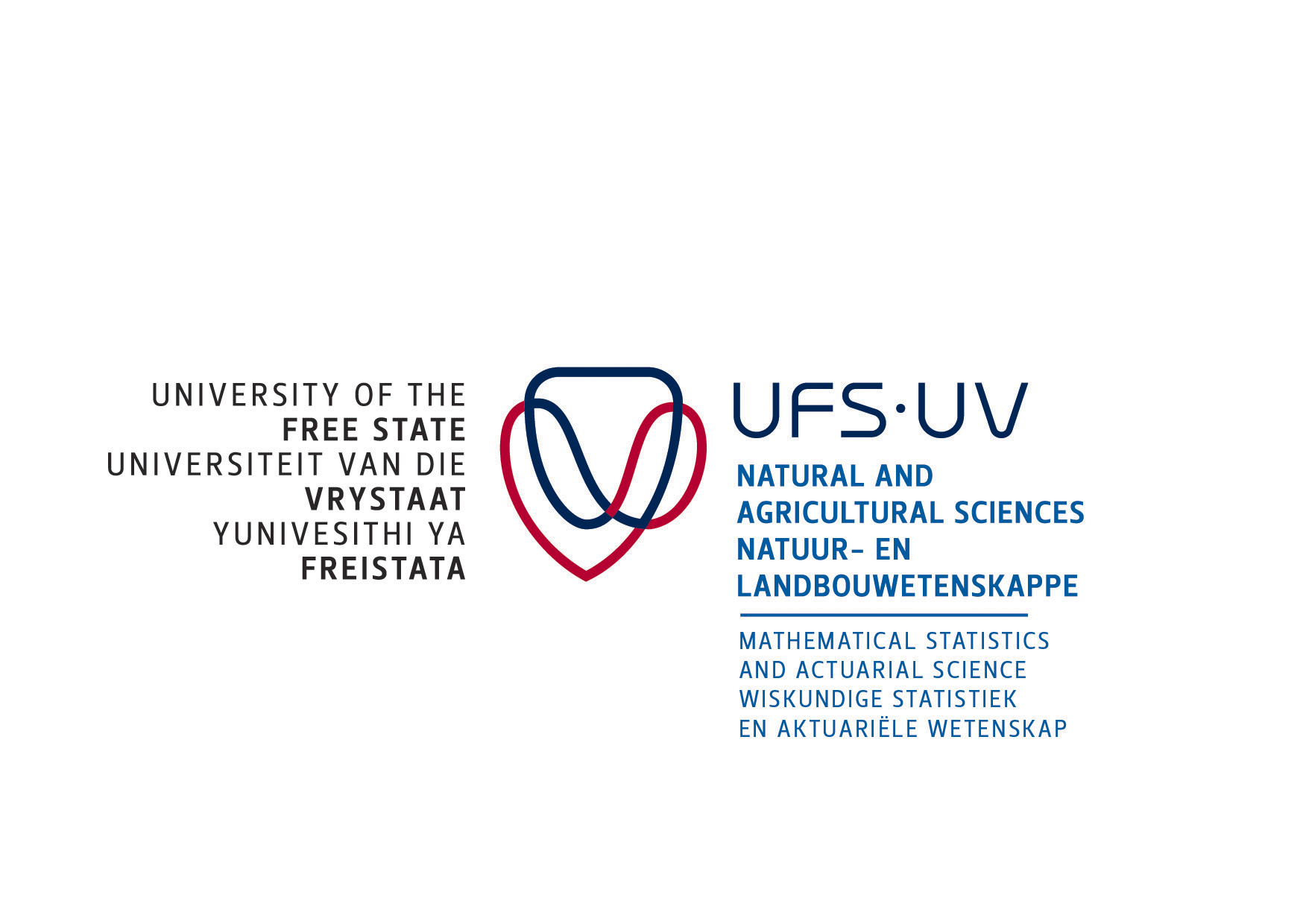}
\end{center}

\textbf{Abstract}

In his 1986 book, Aitchison explains that composition data is regularly mishandled in statistical analyses, a pattern that continues to this day. The Dirichlet Type I distribution is a multivariate distribution commonly used to model a set of proportions that sum to one. Aitchinson goes on to lament the difficulties of Dirichlet modelling and the scarcity of alternatives. While he addresses the second of these issues, we address the first. The Dirichlet Type II distribution is a transformation of the Dirichlet Type I distribution and is a multivariate distribution on the positive real numbers with only one more parameter than the number of dimensions. This property of Dirichlet distributions implies advantages over common alternatives as the number of dimensions increase. While not all data is amenable to Dirichlet modelling, there are many cases where the Dirichlet family is the obvious choice. We describe the Dirichlet distributions and show how to fit them using both frequentist and Bayesian methods (we derive and apply two objective priors). The Beta distribution is discussed as a special case. We report a small simulation study to compare the fitting methods. We derive the conditional distributions and posterior predictive conditional distributions. The flexibility of this distribution family is illustrated via examples, the last of which discusses imputation (using the posterior predictive conditional distributions).


\newpage
\section{Introduction} \label{sec:intro}
The Dirichlet distribution, the multivariate counterpart to the Beta distribution, is often used to model a set of proportions that sum to one, or less than one. In the frequentist sense the distribution is well-documented and researched ---~\citet{Minka}, \citet{Ronning1989}, \citet{Wicker2008}, \citet{Narayanan1992}, \textit{etc}.\ all discuss the various approaches to maximising the likelihood, as well as the method of moments; but Bayesian analysis of this distribution appears to focus entirely on its use as a prior and posterior for the Multinomial distribution. We intend to show that the Dirichlet distribution is useful in its own right and that the Bayesian approach adds to its value.

From \citet{Minka} we see that if $\mathbf{X}=X_1,\ldots,X_{P} \sim Dirichlet(k_1,\ldots,k_{P})$, then the density is
\begin{equation} \label{eq:density}
f(\mathbf{x}) = \frac{\Gamma(k_0)}{\prod_{i=1}^{P}\Gamma(k_i)}\left(\prod_{i=1}^{P}x_i^{k_i-1}\right)
\end{equation}
where $0 < x_i < 1 , i=1,\ldots,P;\  \sum_{i=1}^{P}x_i =1 $ and $ k_0=k_1+\cdots+k_{P}$.

We can say that this distribution is $(P-1)$ dimensional because if any $(P-1)$ values of $\mathbf{x}$ are known then so is the last value. However, the practical relevance and implications of this vary between specific applications.

In the case where $P=2$ the density (\autoref{eq:density}) simplifies to the density of the well-known Beta distribution:
\begin{equation} 
f(x_1) = \frac{\Gamma(k_1+k_2)}{\Gamma(k_1)\Gamma(k_2)}x_1^{k_1-1}(1-x_1)^{k_2-1}
\end{equation}

In general, the marginal distributions are: $X_i\sim Beta(k_i,k_0-k_i) ,i=1,\ldots,P$, with $E(X_i)=\frac{k_i}{k_0}$. A more interesting property (used to derive the above result on the marginal distributions) is the aggregation property \citep[proved in][]{Frigyik2010}. It says that if two components of a Dirichlet vector are dropped and replaced with their sum then the resulting vector is again Dirichlet with the corresponding parameters replaced in the same way. More explicitly,
\begin{equation}\begin{aligned}
&\text{if }(X_1,\ldots,X_P) \sim D(\alpha_1,\ldots,\alpha_P)\\
&\text{then }(X_1,\ldots,X_i+X_j,\ldots,X_P) \sim D(\alpha_1,\ldots,\alpha_i+\alpha_j,\ldots,\alpha_P ).
\end{aligned}\end{equation}
Another way to interpret this property is that the joint marginal distributions are also Dirichlet, in the sense that if $(X_1,\ldots,X_{P-1},1-\sum_{i=1}^{P-1} X_i) \sim D(\alpha_1,\ldots,\alpha_{P-1},\alpha_{P})$, then $(X_1,\ldots,X_j,1-\sum_{i=1}^j X_i) \sim D(\alpha_1,\ldots,\alpha_j,\sum_{i=j+1}^{P} \alpha_i)$. This result becomes very useful later when we discuss the conditional distribution.

The Dirichlet distribution is illustrated in \citet[127--128]{Kotz2000}. Data sets of composition data are given in \citet[354--405]{Aitchison1986} and we use some of these to illustrate the methods in this paper.

In \autoref{sec:fitting} we discuss methods of estimating the parameters of the distribution. In \autoref{sec:type2} we discuss the Dirichlet Type II distribution. This is followed by examples (\autoref{sec:examples}) and a simulation study comparing parameter estimation methods (\autoref{sec:simstudy}). In \autoref{sec:ppcd} we explain the process of prediction using this distribution and in \autoref{sec:imputation} we show how to do imputation using this distribution family.

\section{Fitting the Dirichlet distribution} \label{sec:fitting}

Let $\mathbf{x}_1,\ldots,\mathbf{x}_n$ be an observed sample from the Dirichlet Type I distribution. We combine these vectors as rows of a matrix, so that from here on $X_{ij}$ refers to the $j^{th}$ element of the $i^{th}$ observation. In general, $X_j$ where $j \in {1,\ldots,P}$ will refer to a specific margin considered as a random variable.

\subsection{Method of moments} \label{sec:mom}
\citet{Minka} points out that $k_0=\frac{E(X_1 )-E(X_1^2 )}{E(X_1^2 )-E(X_1 )^2 }$. By multiplying each expected value by this expression we obtain estimates for each parameter. Explicitly:
\begin{equation} 
E(X_j )=\frac{k_j}{k_0} \Rightarrow k_j=E(X_j ) k_0=E(X_j )\left[\frac{E(X_1 )-E(X_1^2 )}{E(X_1^2 )-E(X_1 )^2 }\right]
\end{equation}

This method does not correspond exactly to the commonly used estimates for the Beta case. An alternative approach is given in \citet{Ronning1989} but is not discussed here as the above method is sufficient for our purpose: using the method of moments as a starting point for the method of maximum likelihood.

\subsection{Method of maximum likelihood} \label{sec:ml}

Assuming an independent sample, we obtain the likelihood by multiplying densities. The resulting likelihood is
\begin{equation}\label{eq:D1likelihoodChp2}
Lik(\mathbf{x}_1,\ldots,\mathbf{x}_n|\mathbf{k}) = \left[\frac{\Gamma(k_0)}{\prod_{j=1}^{P}\Gamma(k_j)}\right]^n\prod_{i=1}^n \left(\prod_{j=1}^{P}x_{ij}^{k_{ij}-1}\right)
\end{equation}

\citet{Minka} suggests the following iteration for maximising the likelihood:
\begin{equation} \label{eq:MLinteration}
\psi(k_j^{new} )=\psi(k_0^{old} )+\frac{1}{n} \sum_{i=1}^n \log{X_{ij}}  
\end{equation}
$\psi(x)$ refers to the digamma function, which occurs throughout this text as both $\frac{d}{dx}  \log \Gamma(x)$ and $E[\log G ]$ where $G\sim Gamma(x,1)$. This equivalence between the derivative of the log gamma function and the expectation of the log of a Gamma random variable is useful in derivations that follow. Further, the derivative of the digamma function, denoted $\psi'(x)$, is called the trigamma function and is also used in the coming sections.

A detailed discussion of the methods available to maximise the likelihood is given in \citet{Huang2008}. We found the method above (\autoref{eq:MLinteration}) to converge quickly in all cases. 

\subsection{Posterior distribution} \label{sec:post}

Using a Bayesian approach, the parameters of the Dirichlet distribution can be estimated from the posterior distribution (\textit{e.g.\ }as the mean or mode of the relevant posterior). First, an appropriate prior distribution must be specified.

\subsubsection{MDI Prior} \label{sec:mdi}

We consider first the maximal data information (MDI) prior \citep[41--53]{zellner1997}. The log of the MDI prior is the expected value of the log density. Since
\begin{equation} 
\log{f(\mathbf{x})}=\log {\Gamma(k_0 )}-\sum_{j=1}^{P} \log\Gamma(k_j ) +\sum_{j=1}^{P}\left[(k_j-1) \log {x_j} \right]   
\end{equation}
we have
\begin{equation} \label{eq:mdi}
 \begin{aligned}
 E[\log f(\mathbf{X}) ]&=\log \Gamma(k_0 )-\sum_{j=1}^{P} \log \Gamma(k_j ) +\sum_{j=1}^{P}[(k_j-1)E(\log X_j )] \\ 
&=\log \Gamma(k_0 )-\sum_{j=1}^{P} \log \Gamma(k_j ) +\sum_{j=1}^{P} \left[\left(k_j-1\right)\left(\psi(k_j)-\psi(k_0)\right)\right]. 
\end{aligned}
\end{equation}

The log likelihood is
\begin{equation} \label{eq:llike}
 \begin{aligned}
 \ell &=n \log \Gamma(k_0 )-n \sum_{j=1}^{P} \log \Gamma(k_j )+\sum_{i=1}^n \sum_{j=1}^{P}\left[\left(k_j-1\right) \log x_ij \right] \\ 
&=n \log \Gamma(k_0 )-n \sum_{j=1}^{P} \log \Gamma(k_j ) +\sum_{j=1}^{P} \left[\left(k_j-1\right) \sum_{i=1}^n \log x_{ij}\right]. 
\end{aligned}
\end{equation}

Adding \autoref{eq:mdi} and \autoref{eq:llike} then gives the log posterior:
\begin{equation} \label{eq:lpost}
 \begin{aligned}
 \log \pi(\mathbf{k}|X) &=(n+1) \log \Gamma(k_0 )-(n+1) \sum_{j=1}^{P} \log \Gamma(k_j ) \\
&\qquad+\sum_{j=1}^{P} \left\{\left(k_j-1\right)\left[\left(\psi(k_j)-\psi(k_0)\right)+ \sum_{i=1}^n \log x_{ij} \right]\right\} + c 
\end{aligned}
\end{equation}
where $c$ is an unknown constant.

\citet{Wicker2008} explain that the likelihood is globally concave. The addition of the prior does not appear to change this property, simplifying the process of finding the mode greatly, as the mode is simply the set of parameter values where the posterior reaches its peak. We employed a simple gradient ascent algorithm which only requires the first derivatives:
\begin{equation}
 \begin{aligned}
\frac{\partial \log \pi}{\partial k_j}&=(n+1)[\psi(k_0 )-\psi(k_j )]+\sum_{i=1}^n \log  x_{ij} +(\psi(k_j )-\psi(k_0 ))\\
&\qquad +(k_j-1) \psi' (k_j )-\sum_{i=1}^{P}(k_i-1) \psi' (k_0 )\\
&=n[\psi(k_0 )-\psi(k_j )]+\sum_{i=1}^n \log  x_{ij} +(k_j-1) \psi' (k_j ) - \psi' (k_0 )(k_0-P).
\end{aligned}
\end{equation}

Alternatively, we can simulate a sample from the posterior (\autoref{eq:lpost}), which can be done using the Metropolis-Hastings (MH) algorithm, and calculate the mean. See \citet[267--301]{Robert2004} for an in-depth discussion of this algorithm.

The MH algorithm requires the specification of a suitable jump distribution for choosing candidate values. We chose to use the Multivariate Normal distribution, \textit{i.e.\ }$\mathbf{k}^c \sim N_{P} (\mathbf{k}^i,\Sigma)$ where $\Sigma$ is a diagonal matrix of constants that affect the acceptance rate and required burn-in period. We chose these constants by crudely matching the peak of the posterior to the peak of a Multivariate Normal distribution, one marginal at a time, using the method of percentiles (finding the percentiles using stepwise linear interpolation).

In this way the algorithm is simplified to:

Accept candidate $\mathbf{k}^c$ as a new observation $\mathbf{k}^{i+1}$ from the posterior if and only if

$\log \pi(\mathbf{k}^c |\mathbf{X})-\log \pi(\mathbf{k}^i |\mathbf{X})>\log u$ where $u$ is drawn randomly from the standard Uniform distribution.

It is worth noting that the above method works very well for any finite number of dimensions, including the Beta case $(P=2)$.

\subsubsection{Jeffreys Prior} \label{sec:jeffreys}

An alternative to the MDI prior is the Jeffreys prior \citep{jeffreys1998} which is derived as follows:

Let $g=\log  f(\mathbf{x})=\log \Gamma(k_0 )-\sum_{j=1}^{P} \log \Gamma(k_j )
+\sum_{j=1}^{P} [(k_j-1)  \log  x_j ] $. Then
\begin{equation}
 \begin{aligned}
&\frac{\partial g}{\partial k_j }=\psi(k_0 )-\psi(k_j )+\log x_j ,\\
&\frac{\partial^2 g}{\partial k_j^2 }=\psi'(k_0 )-\psi'(k_j ),\\
&\frac{\partial^2 g}{\partial k_j \partial k_i }=\psi'(k_0 ),
\end{aligned}
\end{equation}
and thus
\begin{equation}
\begin{aligned}
&\text{Jeffreys Prior} \propto \\
&\left|\begin{matrix}
       \psi'(k_1 )-\psi'(k_0 ) & -\psi'(k_0 ) & \cdots & -\psi'(k_0 )         \\
       -\psi'(k_0 ) & \psi'(k_2 )-\psi'(k_0 )  &   & \vdots \\
       \vdots    &  & \ddots & \vdots\\
	-\psi'(k_0) & \cdots & \cdots & \psi'(k_{P} )-\psi'(k_0 )
     \end{matrix}\right|^{\frac{1}{2}}.
\end{aligned}
\end{equation}

In the $Beta$ case $(P=2)$ we get
\begin{equation}
\sqrt{\psi'(k_1 ) \psi'(k_2 )-\psi'(k_0 )\left[\psi'(k_1 )+\psi'(k_2 )\right]}
\end{equation}

When $P=3$ we get:
\begin{equation}
\left\{\psi'(k_1 ) \psi'(k_2 ) \psi'(k_3 ) - \psi'(k_0)\times\left[\psi'(k_1 )\psi'(k_2 ) + \psi'(k_1 )\psi' (k_3 )\psi'(k_2 )\psi'(k_3 )\right]\right\}^{\frac{1}{2}}
\end{equation}

In general the Jeffreys prior is proportional to the square root of
\begin{equation}
 \begin{aligned}
&\prod_{l=1}^{P} \psi'(k_l ) -\psi'(k_0 ) \sum_{l=1}^{P} \prod_{\substack{j=1\\j \ne l}}^{P}\psi' (k_j )\\
&=\prod_{l=1}^{P}\psi' (k_l )\left[1-\psi'(k_0 ) \sum_{l=1}^{P} \frac{1}{\psi'(k_l )}\right]
\end{aligned}
\end{equation}

Thus, the log Jeffreys prior is given by
\begin{equation}
 \begin{aligned}
&0.5 \log \prod_{l=1}^{P}\psi' (k_l )\left[1-\psi'(k_0 ) \sum_{l=1}^{P} \frac{1}{\psi'(k_l )}\right]\\
&=0.5 \sum_{l=1}^{P} \log \psi' (k_l )+0.5 \log\left[1-\psi'(k_0 ) \sum_{l=1}^{P} \frac{1}{\psi'(k_l )}\right]
\end{aligned}
\end{equation}

so that the log posterior is equal to
\begin{equation}
 \begin{aligned}
0.5 \sum_{j=1}^{P} \log \psi' (k_j ) &+0.5 \log \left[ 1 - \psi'(k_0 ) \sum_{j=1}^{P} \frac{1}{\psi'(k_j )} \right] + n \log \Gamma(k_0 )\\
&-n\sum_{j=1}^{P} \log \Gamma(k_j ) + \sum_{j=1}^{P}\left[(k_j-1) \sum_{i=1}^n \log  x_{ij}\right] +c.
\end{aligned}
\end{equation}

The first derivatives (needed to find the mode) are found to be:
\begin{equation}
 \begin{aligned}
\frac{\partial \log  \pi}{\partial k_j } =&n[\psi(k_0 )-\psi(k_j )]+\sum_{i=1}^n \log  x_{ij} + 0.5 \frac{\psi''(k_j )}{\psi'(k_j ) }\\
&-0.5\left[\psi'' (k_0 ) \sum_{l=1}^{P} \frac{1}{\psi'(k_l ) }-\psi' (k_0 ) \psi''(k_j ) (\psi'(k_j ))^{-2}\right]\\
&\times\left[1-\psi' (k_0 ) \sum_{l=1}^{P} \frac{1}{\psi'(k_l ) }\right]^{-1}.
\end{aligned}
\end{equation}

The method of simulation is the same as for the MDI prior (using the Metropolis algorithm).

For an overview of the derivation of objective priors and simpler examples see \citet{Priors}.

\section{The Dirichlet Type II Distribution} \label{sec:type2}

In univariate analysis, the Beta distribution is popular as it is a conjugate prior for the Bernoulli parameter $p$. Should the researcher wish to analyse the odds of success rather than the probability of success, then a conjugate prior is the Beta prime (also called Beta type 2) distribution \parencite{broderick2014}. It is a transform of the Beta distribution \parencite{johnson1970}, according to the following:
\begin{equation}
\text{If }A\sim Beta(a,b)\text{ then }B=\frac{A}{1-A}\sim Beta'(a,b) \text{ and }C=\frac{b}{a}B \sim F(2a,2b).
\end{equation}

A similar transformation can be done in the multivariate case. This transformation was first introduced by \textcite{tiao1965}. Consider again $\mathbf{X} \sim Dirichlet(k_1,\ldots,k_{P} )$. To be clear, $\mathbf{X}$ refers to an unobserved Dirichlet vector and $X_j$ to a specific element of this random vector. If we transform this distribution using the formula 
\begin{equation} \label{eq:type2transform}
 Y_j=\frac{X_j}{X_{P}}   ,j=1,\ldots,(P-1) \ \text{noting that}\ X_{P}=1-\sum_{l=1}^{P-1} X_l 
\end{equation}
then $\mathbf{Y}$ (size $P-1$ by 1) $\sim Dirichlet\ Type\ II(k_1,\ldots,k_{P} )$.

In the above definition the last component is used as the reference component in the transformation as a matter of convenience. While a reference component is required for this transformation, it does not have to be the last component. In practical problems, careful consideration should be given to choosing the reference component. In the absence of natural reference component, a rule of thumb might be to choose the component with the largest minimum, resulting in the most stable transformation.

According to \textcite{tiao1965}, the Dirichlet Type II Distribution also arises as follows: \\
If we have $P\ \chi^2$ random variables with $2k_1,2k_2,\dots,2k_P$ degrees of freedom respectively, say $Q_1,\dots,Q_p$, then the ratios $Y_j=\frac{Q_j}{Q_{P}}   ,j=1,\ldots,(P-1)$ also form a $Dirichlet\ Type\ II(k_1,\ldots,k_{P} )$ distribution.

The density is now
\begin{equation}
 \begin{aligned}
f(\mathbf{y})&=\frac{\Gamma(k_0 )}{\prod_{j=1}^{P}\Gamma(k_j ) } (\prod_{j=1}^{P-1} y_j^{k_j-1} ) (1+y_1+\cdots+y_{P-1} )^{-k_0 }\\
&  \text{where}\ y_j>0 ,i=1,\ldots,(P-1).
\end{aligned}
\end{equation}
and the log density:
\begin{equation}
 \begin{aligned}
\log f(\mathbf{y})&=\log\Gamma(k_0 )-\sum_{j=1}^{P}\log\Gamma(k_j ) - k_0  \log (1+\sum_{j=1}^{P-1} y_j )\\
&+\sum_{j=1}^{P-1} [(k_j-1)  \log y_j] 
\end{aligned}
\end{equation}

Thus, given a sample $(\mathbf{y}_1,\ldots,\mathbf{y}_n )$, the log likelihood is
\begin{equation}
 \begin{aligned}
n \log \Gamma(k_0 )&-n\sum_{j=1}^{P} \log\Gamma(k_j ) -k_0 \sum_{i=1}^n \log \left(1+\sum_{j=1}^{P-1} y_{ij} \right)\\
& +\sum_{i=1}^n \sum_{j=1}^{P-1} [(k_j-1)  \log  y_{ij} ] .
\end{aligned}
\end{equation}

Again we consider the MDI prior. The log prior is:
\begin{equation}
 \begin{aligned}
E[\log f(\mathbf{Y})]&=\log \Gamma(k_0 )-\sum_{j=1}^{P}\log \Gamma(k_j ) - k_0 E\left[\log \left(1+\sum_{j=1}^{P-1} Y_j \right)\right]\\
&+\sum_{j=1}^{P-1} [(k_j-1)E(\log Y_j)] ,
\end{aligned}
\end{equation}
which can be simplified using the following relationships:
\begin{equation}
 \begin{aligned}
E(\log  Y_j )&=E[\log X_j -\log X_P ] \\
& \ \text{ where } X_j \sim Beta(k_j,k_0-k_j) \text{ and } X_P \sim Beta(k_P,k_0-k_P) \\
&=\psi(k_j )-\psi(k_P )
\end{aligned}
\end{equation}
and
\begin{equation}
 \begin{aligned}
&E\left[\log \left(1+\sum_{i=1}^{P-1} Y_i \right) \right] = -E\left[\log \left(1-\sum_{i=1}^{P-1} X_i \right) \right]\\
&=-E[\log X_{P}] \\ 
&=\psi(k_0 )-\psi(k_{P} ).
\end{aligned}
\end{equation}

The log prior thus simplifies to 
\begin{equation}
 \begin{aligned}
\log \Gamma(k_0 )&-\sum_{j=1}^{P} \log \Gamma(k_j ) - k_0 [\psi(k_0 )-\psi(k_{P} )]\\
&+\sum_{j=1}^{P-1} (k_j-1)[\psi(k_j )-\psi(k_P )] 
\end{aligned}
\end{equation}
and the log posterior to
\begin{equation}
 \begin{aligned}
(n+1)&\left[\log \Gamma(k_0 )-\sum_{j=1}^{P} \log \Gamma(k_j ) \right]\\
&-k_0 \left\{\left[\psi(k_0 )-\psi(k_{P} )\right]+\sum_{i=1}^n \log \left(1+\sum_{j=1}^{P-1} y_{ij} \right) \right\}\\
&+\sum_{j=1}^{P-1} (k_j-1)[\psi(k_j )-\psi(k_P)] \\
&+\sum_{i=1}^n \sum_{j=1}^{P-1} [(k_k-1)  \log  y_{ij} ] 
\end{aligned}
\end{equation}

While it is possible to work directly with the above distribution, including simulating from the posterior as before, we prefer to transform a sample from the Dirichlet Type II distribution to a sample from the Dirichlet Type I distribution and apply the methods discussed in \autoref{sec:fitting}. The transformation is the inverse of the previous transformation (\autoref{eq:type2transform}), \textit{i.e.\ } 

\begin{equation}
 X_j=\frac{Y_j}{1+\sum_{l=1}^{P-1} Y_l }  ,j=1,\ldots,(P-1).
\end{equation}

This approach is illustrated in \autoref{sec:examples}.

\section{Examples} \label{sec:examples}

As an example of composition data that should fit the assumptions of the Dirichlet Type I distribution, we consider the hair colours of boys in different counties of Scotland. This data is given in \citet[p. 405]{Aitchison1986}. What identifies this data as Dirichlet is that we expect no additional relationships between hair colours beyond that which is induced by the constraint that the proportions sum to one. This expectation does not hold for the Dark and Black hair colours as these seem to have a strong positive relationship, hence we consider Black hair as part of the Dark hair category prior to applying Dirichlet modelling. After this adjustment the expectation appears to hold.

\begin{figure}
\caption{Hair colour data with simulated replicates}
\begin{center}
\includegraphics[width=\figwidth]{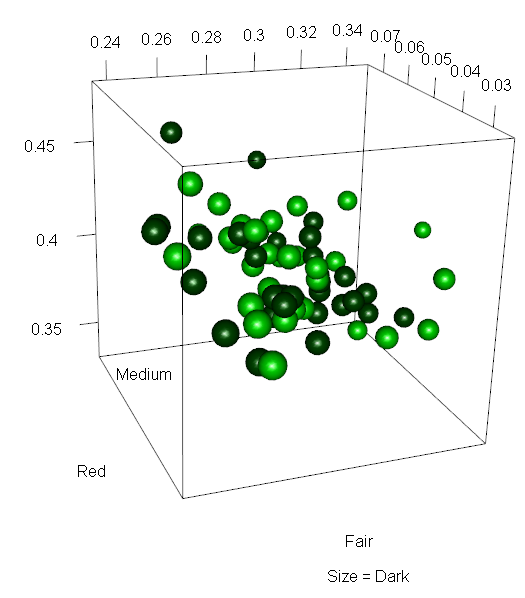}
\end{center}
Dark balls are the observed values and the lighter balls are the simulated values. \label{fig:hair}
\end{figure}

We assess the model fit in three ways. The first is to simulate a replicate sample and compare it to the data visually. This is done in \autoref{fig:hair}. The fit seems reasonable in all 4 dimensions.

The second is to compare the observed and simulated covariance structures (\autoref{tbl:covhair}). These seem remarkably similar, especially when one considers that the Dirichlet has no parameters to explicitly model dependence. The only difference is the apparent positive correspondence between the Fair and Red hair colours, which the Dirichlet model cannot capture.

\begin{table}
\caption{Correlation structure of observed and simulated data}
\begin{center}
\begin{tabular}{|c|c|c|c|c|}
\hline
\textbf{Observed} & Fair  & Red   & Medium & Dark \\
\hline
Fair  & 1     & 0.22  & -0.61 & -0.41 \\
\hline
Red   & 0.22  & 1     & -0.11 & -0.33\\
\hline
Medium & -0.61 & -0.11 & 1     & -0.45 \\
\hline
Dark  & -0.41 & -0.33 & -0.45 & 1 \\
\hline
\textbf{Simulated} & Fair  & Red   & Medium & Dark \\
\hline
Fair  & 1     & -0.14 & -0.52 & -0.39 \\
\hline
Red   & -0.14 & 1     & -0.2  & -0.14 \\
\hline
Medium & -0.52 & -0.2  & 1     & -0.47 \\
\hline
Dark  & -0.39 & -0.14 & -0.47 & 1 \\
\hline
\end{tabular}
   \label{tbl:covhair}
\end{center}
\end{table}

\begin{figure}
\caption{QQ Plots for margins of hair colour data}
\begin{center}
\includegraphics[width=\figwidth]{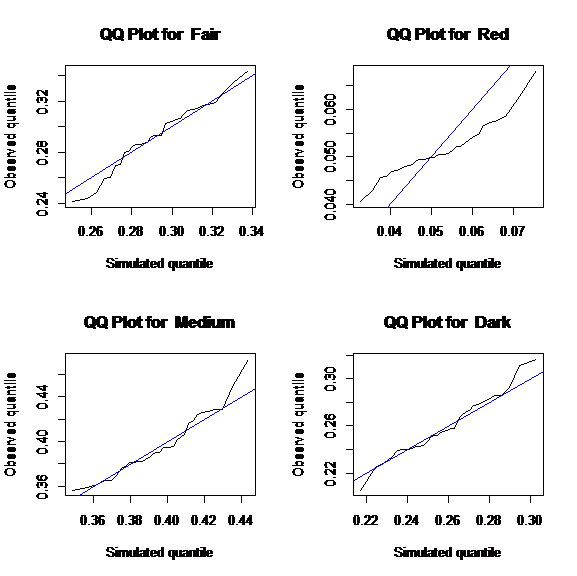}
  \label{fig:qqmargins}
\end{center}
\end{figure}

Lastly, we can view the marginal fits (\autoref{fig:qqmargins}). Again the fit seems reasonable except for the Red hair colour.

As an example of Dirichlet Type II fitting we consider Fisher’s Iris flower data. This is a well-known dataset with 4 variables and three populations of size 50 each. It contains observations on the four variables $Y_1$= sepal length, $Y_2$= sepal width, $Y_3$= petal length and $Y_4$= petal width of the species Setosa, Versicolor and Virginica. This data set is available in the R package \citep{RCore}.

We first assume that the distribution of a set of measurements
\begin{equation}
 \begin{aligned}
&(X_{i1} X_{i2},X_{i3},X_{i4} )=\\
&(\frac{Y_{i1}}{1+\sum_{j=1}^4 Y_{ij} },\frac{Y_{i2}}{1+\sum_{j=1}^4 Y_{ij} },\frac{Y_{i3}}{1+\sum_{j=1}^4 Y_{ij} },\frac{Y_{i4}}{1+\sum_{j=1}^4 Y_{ij} })
\end{aligned}
\end{equation}
can be described through a Dirichlet distribution $D(\alpha_1,\ldots,\alpha_5)$ for each population, say $\mathbf{X}^{(p)}\sim\\ D_4  [\alpha^{(p)}  =(\alpha_j^{(p)},j=1,\ldots,5)],p=1,2,3$.

We proceed to fit Dirichlet distributions on the four variables for each species by transforming as above. We begin by applying the method of moments to obtain initial parameter estimates, as this method is not iterative. We then apply the gradient ascent approach to the posterior distribution using the MDI prior (\autoref{eq:lpost}), and arrive at the following estimates of the Dirichlet parameters:
\begin{equation}
 \begin{aligned}
&\hat{\alpha}^{(SET)}  = {604.71,412.57,175.49,27.78,120.02}\\
&\hat{\alpha}^{(VER)}  = {378.31,176.39,270.87,84.17,64.15}\\
&\hat{\alpha}^{(VIR)}  = {427.61,193.26,360.36,131.62,65.37}.
\end{aligned}
\end{equation}

\begin{figure}[htb]
\caption{Iris data with simulated replicate samples on Dirichlet Type 2 scale}
\begin{center}
\includegraphics[width=0.6\figwidth]{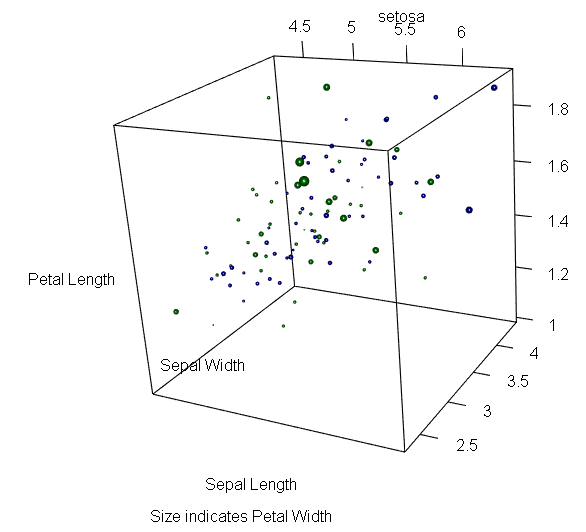}\includegraphics[width=0.6\figwidth]{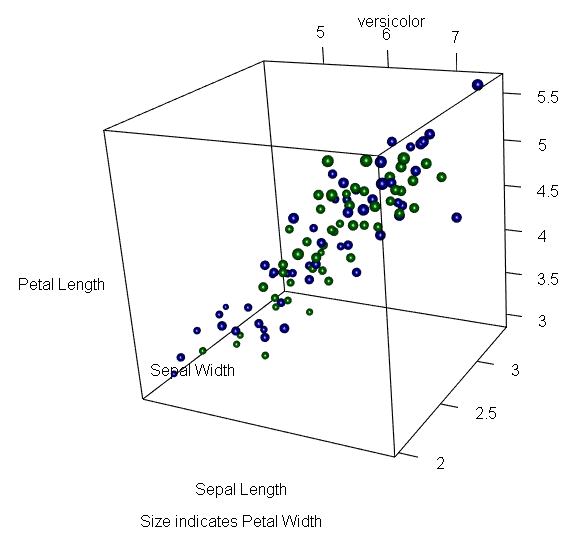}

\includegraphics[width=0.6\figwidth]{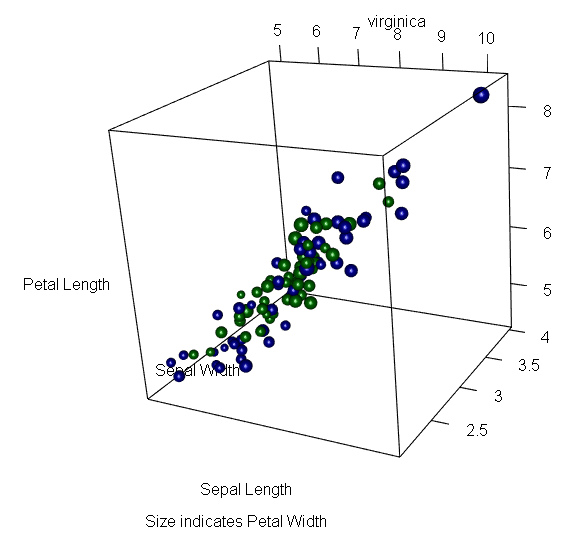}
\end{center}

Green indicates the original data and blue a simulated replicate sample. \label{fig:irisballs}
\end{figure}

A basic test of goodness of fit is to simulate replicate samples and compare them to the original data visually \citep[see for example][159--171]{gelman2013}. Doing this with the iris data produces graphs such as \autoref{fig:irisballs}. There appears to be sufficient correspondence between the observed data and the simulated data. The correspondence is in all four dimensions. More advanced goodness of fit procedures do not yet appear to exist for the Dirichlet distribution.

All methods were implemented in the statistical package R \citep{RCore}. Packages used include the `parallel' package for increasing the speed of the simulation study in \autoref{sec:simstudy} \citep{RCore} and the `rgl' package for drawing the graphs in Figures \ref{fig:hair}, \ref{fig:irisballs} and \ref{fig:bihist} \citep{RGL}.

\section{Simulation Study} \label{sec:simstudy}

We wish to compare the accuracy of the different approaches with respect to estimating the parameters of a Dirichlet distribution. To this end we simulate samples from this distribution for a variety of combinations of sample size, dimension and parameter values. For each combination we simulate 5000 samples. For each sample we estimate the parameters in 6 ways and compare the estimates to the true values. We summarise the discrepancies using the root mean square percentage error (RMSPE) in order to attempt to compare accuracy across varying parameter values.

Simulation of a value from the Dirichlet distribution can be done in at least three ways. The simplest method (used in this study) is to consider a Dirichlet observation as the relative sizes of Gamma draws, \textit{i.e.\ }draw $y_j\sim Gamma(k_j,1)$ for $j=1,\ldots,P$ and then set $x_j=\frac{y_j}{\sum_{j=1}^{P}y_j }$. This gives a draw $ \mathbf{x}\sim D(k_1,\ldots,k_{P} )$.

A second method relies on the fact that the univariate marginal and conditional distributions are Beta and simulates values one dimension at a time conditional on the previous dimensions. A third method by \citet{Matsui2010} uses a path coupling approach. These methods are interesting in the theory of simulation, but not relevant to the current study.

The results of the simulation study are given in \autoref{tbl:simdir}.

\begin{table}
  \caption{Performance of parameter estimation methods as the sample size, parameter magnitude and number of dimensions are varied}
 \footnotesize
\begin{center}
\begin{tabular}{|l|l|l|l|l|l|l|l|}
    \hline
    \textbf{No.} & \textbf{RMSPE} & \textbf{Mode -} & \textbf{ML} & \textbf{Mean -} & \textbf{Mode -} & \textbf{Mean -} & \textbf{MOM} \\
     &  & \textbf{Jeffreys} &  & \textbf{Jeffreys} & \textbf{MDI} & \textbf{MDI} &  \\     \hline
    1     & \textbf{n=25} & 21.90\% & 25.20\% & 25.90\% & 27.20\% & 34.10\% & 38.60\% \\
        & k=(3,3,3) &  &  &  &  &  &  \\    \hline
    2     & \textbf{n=25} & 16.80\% & 19.00\% & 19.70\% & 20.90\% & 26.00\% & 43.00\% \\
        & k=(3,3,3,3,3) &  &  &  &  &  &  \\  \hline
    3     & \textbf{n=50} & 15.10\% & 16.30\% & 16.60\% & 17.00\% & 19.80\% & 23.90\% \\
        & k=(3,3,3) &  &  &  &  &  &  \\ \hline
    4     & \textbf{n=50} & 11.70\% & 12.40\% & 12.70\% & 13.00\% & 15.00\% & 25.70\% \\
        & k=(3,3,3,3,3) &  &  &  &  &  &  \\  \hline
    5     & \textbf{n=90} & 11.10\% & 11.60\% & 11.80\% & 11.90\% & 13.20\% & 17.10\% \\
        & k=(3,3,3) &  &  &  &  &  &  \\  \hline
    6     & \textbf{n=90} & 8.70\% & 9.10\% & 9.20\% & 9.30\% & 10.20\% & 18.30\% \\
        & k=(3,3,3,3,3) &  &  &  &  &  &  \\  \hline
    7     & \textbf{n=120} & 6.50\% & 6.70\% & 6.80\% & 6.90\% & 7.40\% & 16.30\% \\
        & k=(3,3,3,3,3,3,3,3) &  &  &  &  &  &  \\  \hline
    8     & \textbf{n=25} & 21.20\% & 24.50\% & 25.30\% & 26.70\% & 34.20\% & 38.20\% \\
        & k=(4,4,4) &  &  &  &  &  &  \\      \hline
    9     & \textbf{n=25} & 35.30\% & 31.50\% & 30.10\% & 36.00\% & 34.90\% & 38.90\% \\
        & k=(100,100,100) &  &  &  &  &  &  \\  \hline
    10    & \textbf{n=25} & 31.60\% & 27.40\% & 25.70\% & 33.60\% & 34.20\% & 41.60\%\\
        & k=(4,20,100) &  &  &  &  &  &  \\  \hline
    11    & \textbf{n=50} & 15.10\% & 16.30\% & 16.80\% & 17.10\% & 20.00\% & 24.00\% \\
        & k=(4,4,4) &  &  &  &  &  &  \\  \hline
    12    & \textbf{n=50} & 19.30\% & 19.20\% & 18.90\% & 19.90\% & 21.00\% & 23.60\% \\
        & k=(100,100,100) &  &  &  &  &  &  \\  \hline
    13    & \textbf{n=50} & 16.70\% & 17.20\% & 17.10\% & 18.10\% & 20.00\% & 26.10\% \\
        & k=(4,20,100) &  &  &  &  &  &  \\  \hline
    14    & \textbf{n=90} & 10.90\% & 11.40\% & 11.60\% & 11.70\% & 12.90\% & 16.40\% \\
        & k=(4,4,4) &  &  &  &  &  &  \\  \hline
    15    & \textbf{n=90} & 12.00\% & 13.20\% & 13.10\% & 12.40\% & 14.20\% & 16.40\% \\
        & k=(100,100,100) &  &  &  &  &  &  \\  \hline
    16    & \textbf{n=90} & 11.10\% & 12.00\% & 12.10\% & 11.80\% & 13.40\% & 18.40\%\\
        & k=(4,20,100) &  &  &  &  &  &  \\  \hline
    \end{tabular}
\end{center}

\normalsize
\vspace{2mm}
The effect of sample size can be seen by moving from row 1 to 7 or 8 to 16. The effect of parameter magnitude is evident in rows 8 vs 9 vs 10, 11 vs 12 vs 13 or 14 vs 15 vs 16. The effect of the number of dimensions can be seen by comparing rows 1 vs 2, 3 vs 4 or 5 vs 6 vs 7. Lower values indicate better estimation.
 \label{tbl:simdir}
\end{table}

There are many things worth noting in \autoref{tbl:simdir}. The first is that in almost all cases the errors increase as we move from left to right in the table. This means that, in general, the posterior mode using the Jeffreys prior appears to perform the best, followed by the maximum likelihood estimator.

The second is that as some or all parameter values become large then the mean of the posterior simulations (where the mean is taken on each dimension individually) performs better. Thus, if another method is applied and the suggested parameter values are large then this method is recommended.  Special care must be taken when parameter values become small.

The third is that the accuracy of all estimates appears to be heavily dependent on sample size. Increased sample size can be expected to yield increased accuracy but the extent to which it does so is remarkable.

The fourth is that increasing the dimension appears to increase the accuracy of the estimates. Increasing the dimension while keeping the sample size constant increases the total information available from which to estimate parameters. The increase in the information is greater than the increase in parameters, in a relative sense.

We also specifically investigated the performance of these methods for the Beta distribution. We vary the sample size and parameter values. The results are given in \autoref{tbl:simbeta}.

\begin{table}[hbt]
\caption{Performance of parameter estimation methods in the case of the Beta distribution}
\begin{center}
\begin{small}
\begin{tabular}{|l|l|l|l|l|l|l|l|}
    \hline
    \textbf{No.} & \textbf{RMSPE} & \textbf{Mode -} & \textbf{ML} & \textbf{Mean -} & \textbf{Mode -} & \textbf{Mean -} & \textbf{MOM} \\
     &  & \textbf{Jeffreys} &  & \textbf{Jeffreys} & \textbf{MDI} & \textbf{MDI} &  \\      \hline
    1     & \textbf{n=60} & 18.10\% & 19.50\% & 20.40\% & 20.10\% & 23.10\% & 19.60\% \\
    & k=(6,6) &  &  &  &  &  &  \\      \hline
    2     & \textbf{n=60} & 21.10\% & 21.90\% & 24.10\% & 22.40\% & 26.70\% & 21.90\% \\
        & k=(60,60) &  &  &  &  &  &  \\      \hline
    3     & \textbf{n=60} & 20.70\% & 20.70\% & 22.30\% & 20.70\% & 23.10\% & 20.70\% \\
        & k=(600,600) &  &  &  &  &  &  \\      \hline
    4     & \textbf{n=60} & 19.80\% & 21.30\% & 22.80\% & 21.80\% & 26.00\% & 22.60\% \\
        & k=(6,60) &  &  &  &  &  &  \\      \hline
    5     & \textbf{n=25} & 30.50\% & 36.60\% & 38.30\% & 39.00\% & 50.00\% & 36.60\% \\
        & k=(3,3) &  &  &  &  &  &  \\    \hline
    6     & \textbf{n=50} & 20.80\% & 23.00\% & 23.30\% & 23.80\% & 28.40\% & 23.10\% \\
        & k=(3,3) &  &  &  &  &  &  \\    \hline
    7     & \textbf{n=75} & 16.10\% & 17.10\% & 17.50\% & 17.60\% & 20.10\% & 17.20\% \\
        & k=(3,3) &  &  &  &  &  &  \\    \hline
    8     & \textbf{n=100} & 14.30\% & 15.00\% & 15.20\% & 15.30\% & 17.00\% & 15.30\% \\
        & k=(3,3) &  &  &  &  &  &  \\    \hline
    \end{tabular}
\end{small}    
\end{center}

\vspace{1mm}
Parameter magnitude is varied over rows 1 to 4, while sample size is varied over rows 5 to 8.
Lower values indicate better estimation.
  \label{tbl:simbeta}
\end{table}

The patterns are the same as was observed in \autoref{tbl:simdir} --- in all rows we see that the posterior mode using the Jeffreys prior produced the most accurate estimates. 

It should be noted that neither simulation analysis was extended to cover the case where the parameters are smaller than one. In that situation the shape of the posterior distribution may change and force a change in the simulation approach. Parameter values smaller than one correspond with extreme variance, with observations more spread out than one would expect from a Uniform distribution. Extra care must be taken if such data is encountered in practice.

\section{Predictive posterior conditional distribution} \label{sec:ppcd}

We consider the situation where a number of components are jointly distributed Dirichlet and we wish to predict the values of a subset of components, given values for the remainder of the components. We can do this using the predictive posterior conditional distribution.

\subsection{Conditional distribution}
\begin{sloppypar}Initially assume that the parameter values are known, \textit{i.e.\ } $(X_1,\ldots,X_P )|k \sim D(k_1,\ldots,k_{P} )$ and we are interested in the distribution of $\left(X_1,\ldots,X_{P-1} \middle| X_{P} \right)$. We condition on the last component for convenience only, the result holds for any component. We use the result that the last component follows a $Beta(k_P,\sum_{j=1}^{P-1})$ marginal distribution, as proved in \textcite{Frigyik2010}.\end{sloppypar}

\begin{equation}
 \begin{aligned}
f(x_1,\ldots,x_{P-1}|X_P=x_P )&=\frac{f(x_1,\ldots,x_P )}{f(x_P)} \\
&=\frac{\frac{\Gamma\left(\sum_{j=1}^{P} k_j \right)}{\prod_{j=1}^{P} \Gamma(k_j ) } \left(\prod_{j=1}^P x_j^{k_j-1} \right)}{\frac{\Gamma\left(\sum_{j=1}^{P} k_j\right) }{\Gamma(k_P ) \Gamma\left(\sum_{j=1}^{P-1} k_j\right)} \left( x_P^{k_P-1} (1-x_P)^{\left(\sum_{k=1}^{P-1} k_j\right) -1} \right)}\\
&=\frac{\Gamma\left(\sum_{j=1}^{P-1} k_j\right)}{\prod_{j=1}^{P-1} \Gamma(k_j ) } \left(\prod_{j=1}^{P-1} \left(\frac{x_j}{1-x_P}\right)^{k_j-1} \right) (1-x_P)^{1-P}
\end{aligned}
\end{equation}

If we set $W_j=\frac{X_j}{1-x_P}$, where $j=1,\dots,(P-1)$ and $x_P$ is a constant in this situation, then $f(w_1,\ldots,w_{P-1}|X_{P}=x_{P} )$ becomes
\begin{equation}
\frac{\Gamma\left(\sum_{j=1}^{P-1} k_j\right)}{\prod_{j=1}^{P-1} \Gamma(k_j ) } \left(\prod_{j=1}^{P-1} \left(\frac{x_j}{1-x_P}\right)^{k_j-1} \right)
\end{equation}
since the absolute Jacobian is $(1-x_P)^{P-1}$.

\begin{sloppypar}From this expression we see that $(W_1,\ldots,W_{P-1} |X_{P}=x_{P} ) \sim D(k_1,\ldots,k_{P-1} )$.\end{sloppypar}

The result above can be applied recursively, each time eliminating one component by dropping its parameter and rescaling the remaining components by their sum to ensure they again sum to one, as required by the Dirichlet distribution.

Thus, given a set of parameters and a set of known values for some components we can easily calculate quantities of interest. For example, we can simulate from the conditional distribution to calculate prediction intervals.

\subsection{Implementation of the predictive posterior conditional distribution}
In practise the parameters are unknown and must be estimated from data. This carries uncertainty across to the predictions, uncertainty that must be taken into account. We adjust for the uncertainty using the predictive posterior conditional distribution.

Practically, the process for a single vector observation works as follows:
\begin{enumerate}
 \item Given a complete data matrix $X^{(c)}$ of size $n$ by $P$, we simulate a large number of sets of parameters from the posterior distribution. Suppose we simulate $T$ sets and store them in a $T$ by $P$ parameter matrix $K$.
 \item Now consider a new matrix $X^{(n)}$ of size $m$ by $P$, where $Q$ components require prediction $(1\leq Q <P)$. For convenience we can arrange the components of the new matrix, as well as the simulated parameter matrix, such that the $Q$ components are grouped together on the left.
 \item For each individual vector of simulated parameters $(\mathbf{k}_t)$, we simulate a set of $w_1,\dots,w_{Q}$ from a $D(K_{t1},\ldots,K_{tQ})$, where $t=1,\ldots,T$.
 \item Multiply each simulated set of $w_1,\dots,w_{Q}$ by $\prod_{j=Q+1}^P (1-X^{(n)}_j)$.
 \item Combine the predictions over simulations  $t=1,\ldots,T$ to calculate values of interest, since these simulations constitute an approximation of the predictive posterior conditional distribution.
\end{enumerate}

\subsection{Conditional distribution in the Type II case}

In the case where the given values are on the Dirichlet Type II scale we follow a different approach where we apply the aggregation property directly:

We know that $ \left(X_1,\ldots,X_j,\left[1-\sum_{l=1}^j X_l \right]\right)  \sim D\left(k_1,\ldots,k_j,\sum_{l=j+1}^{P} k_l \right)$ and that $X_j=\frac{Y_j}{1+\sum_{l=1}^{P-1} Y_l }$, where $\left(Y_1,\ldots,Y_{P-1}\right) \sim D2(k_1,\ldots,k_{P})$.

Thus, given $(Y_{j+1},\ldots,Y_{P-1} )$, we have that
\begin{equation}
\frac{\frac{(Y_1,\ldots,Y_j )}{1+\sum_{l=1}^{P-1} Y_l }}{1-\sum_{i=1}^j \frac{Y_i}{1+\sum_{l=1}^{P-1} Y_l } } \sim D2\left(k_1,\ldots,k_j,\sum_{l=j+1}^{P} k_l \right).	
\end{equation}

This expression simplifies to
\begin{equation}
\frac{(Y_1,\ldots,Y_j )}{1+\sum_{l=j+1}^{P-1} Y_l }  ,
\end{equation}
so, if we define $\beta=\left(1+\sum_{l=j+1}^{P-1} Y_l \right)$ and $ V_i=\frac{Y_i}{\beta} $then
\begin{equation}
(V_1,\ldots,V_j |Y_{j+1}=y_{j+1},\ldots,Y_{P-1} = y_{P-1}) \sim D2\left(k_1,\ldots,k_j,\sum_{l=j+1}^{P} k_l \right).
\end{equation}

We can also define $U_i=\frac{V_i}{1+\sum_{l=1}^j V_l }$ and note that
\begin{equation}
\left(U_1,\ldots,U_j,1-\sum_{l=1}^j U_l \middle|Y_{j+1}=y_{j+1},\ldots,Y_{P-1}=y_{P-1}\right) \sim D\left(k_1,\ldots,k_j,\sum_{l=j+1}^{P} k_l \right).
\end{equation}

We can use either of these and appropriately modify steps 2 and 3 in the given process to calculate desired quantities. For example, if we simulate a set of $u_1,\dots,u_j$ then we note that $y_l=\frac{\beta u_l}{1-\sum_{h=1}^j u_h }$.

\begin{figure}[hbt]
\caption{Bivariate histogram to illustrate the predictive posterior conditional distribution}
\begin{center}
\includegraphics[width=\figwidth]{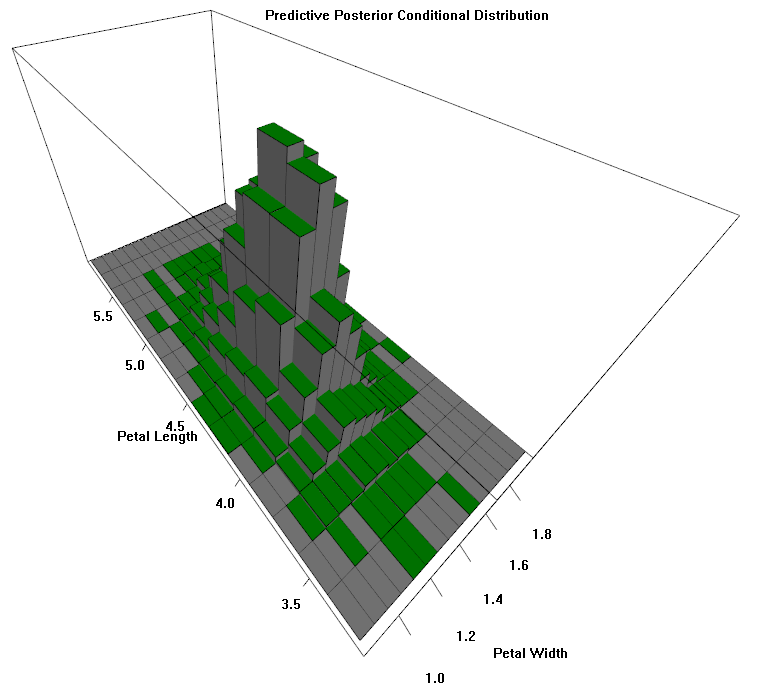}
\end{center}

Distribution of the petal size given that the sepal length and width are 6 and 3 respectively, for the Versicolor species. \label{fig:bihist}
\end{figure}

As an example, let us again consider the Iris data from \autoref{sec:examples}. Suppose we observe a new flower from the Versicolor family with a sepal length of 6 and a sepal width of 3, but we are unable to measure the petals. We apply the above process with $T=9999$ and find that the predictive posterior conditional distribution has the form depicted in \autoref{fig:bihist}. By taking means we predict a petal length of 4.38 and a petal width of 1.36, which seems reasonable looking at \autoref{fig:irisballs}.

\section{Multiple Imputation of Dirichlet Type II Data} \label{sec:imputation}

As an additional example we include an application of the methods discussed. We consider a set of monthly total rainfall figures collected at five stations in central South Africa. We only consider months where measurable amounts of rain fell. We have 37 months where all five stations produced measurements and 58 months where a station’s data is missing and needs to be imputed.

While many imputations methods exist, it is generally agreed that a flexible and efficient solution is multiple imputation. For an overview of multiple imputation see \citet{Schafer1998}. 

We assume that the rainfall measurements follow a Dirichlet Type II distribution with unknown parameters. We proceed through the following steps:

\begin{enumerate}
  \item Estimate the six parameters of the distribution using only the 37 complete months.
  \item Pick random values for all the missing measurements from the respective predictive posterior conditional distribution of each. This is done one at a time since each distribution is different.
 \item Estimate the six parameters of the distribution using all 95 months.
 \item Pick random values for all the missing measurements from the respective predictive posterior conditional distribution of each.
 \item Repeat steps 3 and 4 a number of times, preferably until the parameter estimates converge in distribution.
 \item Store the last set of completed data and last set of parameter estimates.
  \item Repeat steps 1 to 6 a number of times (min. 5) to obtain a collection of imputed data sets.
\end{enumerate}

The parameter estimates that are produced as a result of this approach are given in \autoref{tbl:imputes}. It is clear that parameters $k_2, k_4$ and $k_5$ would be underestimated if only the complete data were used.

An extract of imputed data values is given in \autoref{tbl:imputespred} for illustration.

It is clear that this approach provides a valid alternative to the commonly used Multivariate Normal imputation \citep[see][for an explanation of why this is important]{vonHippel2013}. Another, perhaps even more important use of this method is that it can address the problem of zeros in composition data. Zeros can be made missing and then multiply imputed.

\begin{table}[bh]
    \caption{Parameter estimates following multiple imputation} 
\begin{center}
    \begin{tabular}{|l|r|r|r|r|r|r|}
    \hline
    \textbf{Parameter estimates} & \multicolumn{1}{c|}{\textbf{k1}} & \multicolumn{1}{c|}{\textbf{k2}} & \multicolumn{1}{c|}{\textbf{k3}} & \multicolumn{1}{c|}{\textbf{k4}} & \multicolumn{1}{c|}{\textbf{k5}} & \multicolumn{1}{c|}{\textbf{k6}} \\
    \hline
    Complete data only & 3.57  & 3.49  & 2.72  & 2.8   & 3.62  & 0.42 \\    \hline
    Imputation 1 & 3.43  & 3.83  & 2.78  & 3.13  & 3.97  & 0.43\\    \hline
    Imputation 2 & 3.43  & 3.89  & 2.65  & 2.91  & 4     & 0.43 \\    \hline
    Imputation 3 & 3.5   & 3.82  & 2.77  & 3.03  & 4.03  & 0.43 \\    \hline
    Imputation 4 & 3.46  & 3.83  & 2.75  & 3.12  & 4.07  & 0.43 \\    \hline
    Imputation 5 & 3.55  & 3.94  & 2.8   & 3.17  & 4.05  & 0.43 \\    \hline
    Imputation 6 & 3.46  & 3.81  & 2.78  & 3.15  & 3.98  & 0.43 \\    \hline
    Imputation 7 & 3.63  & 4.04  & 2.84  & 3.33  & 4.17  & 0.44 \\    \hline
    Imputation 8 & 3.2   & 3.49  & 2.54  & 2.83  & 3.7   & 0.42 \\    \hline
    Imputation 9 & 3.28  & 3.65  & 2.65  & 2.93  & 3.66  & 0.42 \\    \hline
    Imputation 10 & 3.59  & 4.02  & 2.85  & 2.98  & 4.12  & 0.44 \\    \hline
    \end{tabular}
\label{tbl:imputes}
\end{center}
\end{table}

\begin{table}[thb]
  \caption{Extract from final round of Imputation 2 showing typical imputed values}
\begin{center}
    \begin{tabular}{|l|r|r|r|r|r|}
    \hline
    \textbf{Month} & \multicolumn{1}{l|}{\textbf{Station 1}} & \multicolumn{1}{l|}{\textbf{Station 2}} & \multicolumn{1}{l|}{\textbf{Station 3}} & \multicolumn{1}{l|}{\textbf{Station 4}} & \multicolumn{1}{l|}{\textbf{Station 5}} \\    \hline
    2003m4 & 5.8   & 6.5   & 8.5   & \textbf{6.4} & 9.4 \\    \hline
    2003m8 & 6.3   & 3.4   & 2.1   & \textbf{2.6} & 4.6 \\    \hline
    2003m9 & 30.4  & 18    & 20.2  & \textbf{28.3} & 20.3 \\    \hline
    2003m11 & 47.6  & 49.7  & 16.3  & \textbf{55.1} & 42.8 \\    \hline
    2003m12 & 24.3  & 52.4  & 0.7   & \textbf{7.9} & 45.2 \\    \hline
    2004m1 & 5.6   & 60.6  & 3.2   & \textbf{2.4} & 49.2 \\    \hline
    2004m3 & 73.8  & 85.3  & \textbf{75.7} & 46.7  & 151.5 \\    \hline
    \end{tabular}
   \label{tbl:imputespred}
\end{center}
\end{table}

\section{Conclusion} \label{sec:conclusion}
We have found that the Dirichlet distribution family is a useful distribution when modelling random variables that take on positive values and especially for composition data with negative correlations, particularly as the number of dimensions increases. We see that it scales easily from the univariate (Beta) case to the multivariate (Dirichlet) case with the number of parameters only increasing by one per dimension. We found that the Bayesian approach can add accuracy and flexibility when working with this distribution family. We have looked at conditioning within this distribution family and how this can be extended to allow imputation, which in turn can address the problem of zeros in observed data.

\printbibliography[title=References]

\end{document}